\documentclass[12pt]{elsart}
\usepackage{latexsym}
\usepackage{amssymb}
\usepackage{amsmath}
\usepackage{xypic}
\usepackage{graphicx}

\newcommand{\cB}{\mathcal{B}}
\newcommand{\cC}{\mathcal{C}}
\newcommand{\clD}{\mathcal{D}}
\newcommand{\cE}{\mathcal{E}}

\newcommand{\cN}{\mathcal{N}}

\newcommand{\bN}{\mathbb{N}}

\newcommand{\bW}{\mathbb{W}}
\newcommand{\bZ}{\mathbb{Z}}
\newcommand{\bfa}{{\bf a}}
\newcommand{\bfb}{{\bf b}}
\newcommand{\bfc}{{\bf c}}
\newcommand{\bff}{{\bf f}}

\newcommand{\ass}{\mathfrak{a}}
\newcommand{\gss}{\mathfrak{g}}
\newcommand{\qss}{\mathfrak{q}}
\newcommand{\com}{\mathfrak{c}}
\newcommand{\lid}{\mathfrak{l}}

\newcommand{\rid}{\mathfrak{r}}

\newcommand{\ph}[1]{\phantom{#1}}

\newcommand{\can}{\textrm{{\bf can}}}

\newcommand{\End}{\textrm{End}}

\newcommand{\maxi}{\textrm{max}}
\newcommand{\mini}{\textrm{min}}

\newcommand{\ubar}{\underline{\ph{S}}}
\newcommand{\Funct}{\textrm{\bf Funct}}
\newcommand{\LO}{\textrm{\bf LO}}
\newcommand{\xB}{\textrm{\bf xB}}
\newcommand{\ev}{\textrm{ev}}
\newcommand{\Cptr}{\textrm{\bf Cptr}}
\newcommand{\BCptr}{\textrm{\bf BCptr}}
\newcommand{\NCptr}{\textrm{\bf NCptr}}
\newcommand{\BNCptr}{\textrm{\bf BNCptr}}
\newcommand{\Cohr}{\textrm{\bf Cohr}}

\begin{document}

\title{Braided Premonoidal Coherence}
\author{William P. Joyce}
\address{Department of Physics \& Astronomy, University of Canterbury,
  Private Bag 4800, Christchurch, New Zealand.}

\maketitle

\begin{abstract}
Given a category with a bifunctor and natural isomorphisms for
associativity, commutativity and left and right identity we do not
assume that extra constraining diagrams hold. We introduce groupoids of
coupling trees to describe a version of coherence that is weaker than
the usual notion of monoidal coherence.
\end{abstract}

\section{Introduction}                                               %

The constraints of a monoidal category are too strong for some
applications in physics. Although the natural isomorphism property is
founded on strong physical principles (see Joyce \cite{wj2}) there are
no conclusive grounds for requiring that the pentagon and triangle
diagrams hold. For example in Joyce \cite{wj4,wj3} recouplings
embodying the Pauli exclusion principle are natural isomorphisms that do
not satisfy the pentagon diagram for representations of the group
$SU(n)$ where $n>2$. This provides motivation for investigating weakened
monoidal structures. A partial exploration of the problem is given in
Joyce \cite{wj} and does not consider braids, nor violations of the
triangle diagram. An alternative approach using $n$--categories is given
by Yanofsky \cite{ny}

Given a category $\cC $ with a bifunctor $\otimes :\cC \times \cC
\rightarrow \cC $ one may form ``words'' by iterating the bifunctor.
Thus one may ask what is the relationship between different words. We
initially investigate structures with only a natural isomorphism $\ass
:\otimes (\otimes \times 1)\rightarrow \otimes (1\times \otimes )$ for
associativity. Any extra conditions such as the usual pentagonal
constraint of a monoidal category is a luxury and not necessarily a
fundamental requirement. This minimum structure we call a premonoidal
structure. We introduce the groupoid of coupling trees which
characterises this structure functorially. Further we consider weaker
constraints than the pentagonal constraint and prove suitable coherence
theorems. Ultimately we arrive at the Mac Lane coherence theorem
\cite{sm2} without a unit.

Next we consider two extensions of these structures. First to a category
with a unit structure where the triangle diagram is not necessarily
required to hold. The relevant groupoid is of coupling trees with
nodules. The second extension is to braids where the hexagonal diagrams
are not necessarily required to hold. This requires the exploded
groupoid of braids. The case of braided monoidal coherence is elegantly
addressed by Joyal and Street \cite{ajrs}. The present paper concludes
by exhibiting a diagram calculus representing all commutative diagrams
in $\cC $.

We collect some elementary, but useful results on functor categories. As
a general reference we recommend Mac Lane's book \cite{sm}.
\begin{prop}
Let $\tau :F\rightarrow G$ be a natural transformation between two
functors $F,G:\cC \rightarrow \clD $.
\begin{enumerate}
\item{Given a functor $H:\cB \rightarrow \cC $ then $\tau
    H:FH\rightarrow GH $ defined by $b\mapsto \tau_{Ha}:FHb\rightarrow
    GHb$ is a natural transformation. Furthermore, if $\tau $ is a natural
    isomorphism, then so is $\tau H$.}
\item{Given a functor $K:\clD \rightarrow \cE $ then $K\tau
    :KF\rightarrow KG$ defined by $c\mapsto K\tau_{c}:KFc\rightarrow
    KGc$ is a natural transformation. Furthermore, if $\tau $ is a
    natural isomorphism, then so is $K\tau $.}
\item{Given another natural transformation $\sigma :H\rightarrow K$ then
    $\tau \times \sigma :F\times H\rightarrow G\times K$ is a natural
    transformation defined by $(c,b)\mapsto (\tau_{c},\sigma_{b})$.}
\end{enumerate}
\end{prop}
The proof is straightforward and given by checking the naturality
condition holds using the functorial properties of $H$ and $K$.

\section{Premonoidal Structure of a Category}                        %

A category $\cC $ with a functor $\otimes :\cC \times \cC \rightarrow
\cC $ does not have sufficient properties to be of any interest.
However, including a natural isomorphism for associativity leads to a
surprisingly rich structure. Moreover, this section and the next is a
warm up for the richer structures of later sections.
\begin{defn}
  A premonoidal structure for a category $\cC $ is a doublet $(\otimes
  ,\ass )$ where $\otimes :\cC \times \cC \rightarrow \cC $ is a functor
  called tensor product and $\ass :\otimes (\otimes \times 1)\rightarrow
  \otimes (1\times \otimes )$ is a natural isomorphism called the
  recoupling for associativity.
\end{defn}
A premonoidal structure is not necessarily monoidal because the pentagonal
constraint does not hold and there is no account taken of a unit.
Despite the lack of a pentagonal constraint we can nevertheless measure
the degree to which the pentagonal constraint is deformed as defined in
the next definition.
\begin{defn}
  The recoupling for deformativity is the natural automorphism
  $\qss :\otimes (\otimes \times \otimes )\rightarrow \otimes (\otimes
  \times \otimes )$ defined by
\begin{eqnarray}
\qss_{a,b,c,d} & = & \ass^{-1}_{a,b,c\otimes
  d}(1_{a}\otimes \ass_{b,c,d})\ass_{a,b\otimes c,d}(\ass_{a,b,c}\otimes 
1_{d})\ass^{-1}_{a\otimes b,c,d}
\end{eqnarray}
  for all objects $a,b,c,d$ of $\cC $, as depicted in the following
  diagram.
\end{defn}
\def\objectstyle{\scriptstyle}
\xymatrix@R=60pt@C=25pt{ & (a\otimes b)\otimes (c\otimes d)
  \ar[r]^*+{\qss_{a,b,c,d}} &
  (a\otimes b)\otimes (c\otimes d) \ar[dr]^*+{\ass_{a,b,c\otimes d}} & \\
  ((a\otimes b)\otimes c)\otimes d \ar[ur]^*+{\ass_{a\otimes b,c,d}}
  \ar[dr]_*+{\ass_{a,b,c}\otimes 1_{d}} &
  & & a\otimes (b\otimes (c\otimes d)) \\
  & (a\otimes (b\otimes c))\otimes d \ar[r]_*+{\ass_{a,b\otimes c,d}} &
  a\otimes ((b\otimes c)\otimes d) \ar[ur]_*+{1_{a}\otimes \ass_{b,c,d}} &
  \\ }
Note that if $\qss_{a,b,c,d}=1_{(a\otimes b)\otimes (c\otimes d)}$ for
all objects $a,b,c,d$ then the pentagonal constraint holds. This is then a
monoidal structure without unit. Next we define the pseudo--monoidal
structure first introduced in Joyce \cite{wj}
\begin{defn}
A pseudo--monoidal structure for a category $\cC $ is a premonoidal
structure $(\otimes ,\ass )$ such that the following two dodecagon
diagrams commute.\\ \\
\xymatrix@R=40pt@C=40pt{
((a\otimes b)\otimes c)\otimes (d\otimes f) \ar[r]^{\qss_{a\otimes
    b,c,d,f}} \ar[d]_{\ass_{a,b,c}\otimes 1_{d\otimes f}} & ((a\otimes
  b)\otimes c)\otimes (d\otimes f) \ar[d]^{\ass_{a,b,c}\otimes
    1_{d\otimes f}}\\
(a\otimes (b\otimes c))\otimes (d\otimes f) \ar[r]_{\qss_{a,b\otimes
    c,d,f}} & (a\otimes (b\otimes c))\otimes (d\otimes f)\\}\\ \\ \\
\xymatrix@R=40pt@C=40pt{
(a\otimes b)\otimes ((c\otimes d)\otimes f) \ar[r]^{\qss_{a,b,c\otimes
    d,f}} \ar[d]_{1_{a\otimes b}\otimes \ass_{c,d,f}} & (a\otimes
b)\otimes ((c\otimes d)\otimes f) \ar[d]^{1_{a\otimes b}\otimes
  \ass_{c,d,f}}\\
(a\otimes b)\otimes (c\otimes (d\otimes f)) \ar[r]_{\qss_{a,b,c,d\otimes 
    f}} & (a\otimes b)\otimes (c\otimes (d\otimes f))\\}\\ \\
for all objects $a,b,c,d,f$ of $\cC $.
\end{defn}

\section{Coherence of Premonoidal Structures}                        %

We begin with some preliminary definitions. Let $[n]=\{ 1,2,...,n\} $.
\begin{defn}
  A coupling tree $t$ of length $n$ is a planar binary rooted tree with
  $n$ leaves, together with a linear ordering of its vertices subject to
  the condition that any connected loop--free sequence of vertices from
  the root to a leaf is (strictly) increasing. Hence all but the null
  coupling tree, are uniquely characterised by a bijection
  $t:[n-1]\rightarrow [n-1]$ giving the order in which the branch point
  levels occur. The length of the tree, denoted $|t|$, is the number of
  its leaves $n$.
\end{defn}
Note that the null coupling tree is represented formally by
$0:[-1]\rightarrow [-1]$, where as the empty map $1:[0]\rightarrow [0]$
represents the (unique) coupling tree of length one. We denote the
groupoid of coupling trees of length $n$ by $\Cptr_{n}$ where there is a
unique arrow between each tree called a recoupling. Thus the groupoid of
coupling trees is given by
\begin{eqnarray}
\Cptr & = & \coprod_{n\in \bW }\Cptr_{n}
\end{eqnarray}
where we define $\bW =\bN \cup \{ 0\} $. A coupling tree of length $n$
is equally well represented by a linear ordering of the elements $[n-1]$
in the following obvious way. Writing the level at the bottom of each
region between adjacent leaves from left to right in a sequence gives a
linear ordering. For example, the coupling tree
$(1243)(5):[5]\rightarrow [5]$
of length $6$ has the linear ordering\\ \\
\hspace*{2cm}\includegraphics[width=100pt]{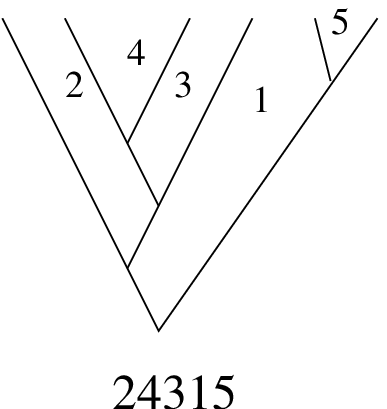}\\ \\
This defines an injective functor from $\LO $ into $\Cptr $ where $\LO $
denotes the groupoid of linear orderings and
\begin{eqnarray}
\LO & = & \coprod_{n\in \bW }\LO_{n}
\end{eqnarray}
where $\LO_{n}$ is the full subgroupoid generated by the linear
orderings of length $n$. The functor between $\LO $ and $\Cptr $ is
invertible if we extend $\LO $ to $\LO \cup \{ *\} $ where $*$ is a
discrete object mapping to the null tree. The recouplings between two
coupling trees (of the same length) are represented by permutations. In
what follows we do not distinguish between the two groupoids using the
linear ordering to denote the coupling trees and recouplings.

If we cut a coupling tree $t$ at its root then we obtain two coupling
trees. We denote the left coupling tree by $Lt$ and the right by $Rt$.
We may also split recouplings (or permutations) about the root vertex.
Let $\pi :s\rightarrow t$ be an arrow of $\Cptr $ that leaves the root
fixed. Furthermore, if $r$ is the position of the root in $s$, $n=|s|-1$
and $[r-1]$ and $r+[n-r]\equiv \{ r+1,...,n\} $ are closed under $\pi $,
then we define $L\pi :Ls \rightarrow Lt$ to be the unique recoupling
from $Ls$ to $Lt$, similarly for $R\pi :Rs\rightarrow Rt$.

We now have the following coherence result where we use the notation
\begin{eqnarray}
\cC^{S} & \equiv & \coprod_{s\in S}\cC^{s}
\end{eqnarray}
for any $S\subset \bW $ with $\cC^{0}$ defined to be the one arrow category.
\begin{thm} \label{cohthm1}
Given a category $\cC $ with premonoidal structure $(\otimes ,\ass )$
there is a unique functor $\Gamma :\Cptr \rightarrow \Funct
(\cC^{\bW },\cC )$ satisfying:
\begin{enumerate}
\item{\begin{eqnarray}
\Gamma (t) & = & \otimes \left( \Gamma (Lt)\times \Gamma (Rt)\right)
\label{propty1} \\
\Gamma (\emptyset ) & = & 1
\end{eqnarray}
for all objects $t$ of $\Cptr $.}
\item{Let $(ij):s\rightarrow t$ be a transposition interchanging $k$ in
the $i$th position with $k+1$ in the $j$th position of $s$ such that
the position of level $1$ is not between $i$ and $j$, then for $k=1$
\begin{eqnarray}
\Gamma (ij) & = & \left\{ \begin{array}{cc} \ass \left( \Gamma (LLs)\times
    \Gamma (RLs)\times \Gamma (Rs) \right) & : i>j\\ \ass^{-1} \left(
    \Gamma (Ls) \times \Gamma (LRs) \times \Gamma (RRs) \right) & : i<j
\end{array} \right.
\end{eqnarray}
for $k>1$ and $s^{-1}(1)<\mini \{ i,j\}$
\begin{eqnarray}
\Gamma (ij) & = & \otimes (\Gamma L(ij) \times 1_{\Gamma (Rs)})
\end{eqnarray}
and for $k>1$ and $s^{-1}>\maxi \{ i,j\}$
\begin{eqnarray}
\Gamma (ij) & = & \otimes (1_{\Gamma (Ls)}\times \Gamma R(ij))
\end{eqnarray}}
\end{enumerate}
\end{thm}
Before giving a proof we make some remarks and introduce some convenient 
notation. The deformed pentagon diagram in $\Cptr $ is given by\\ \\
\hspace*{2cm}\includegraphics[width=240pt]{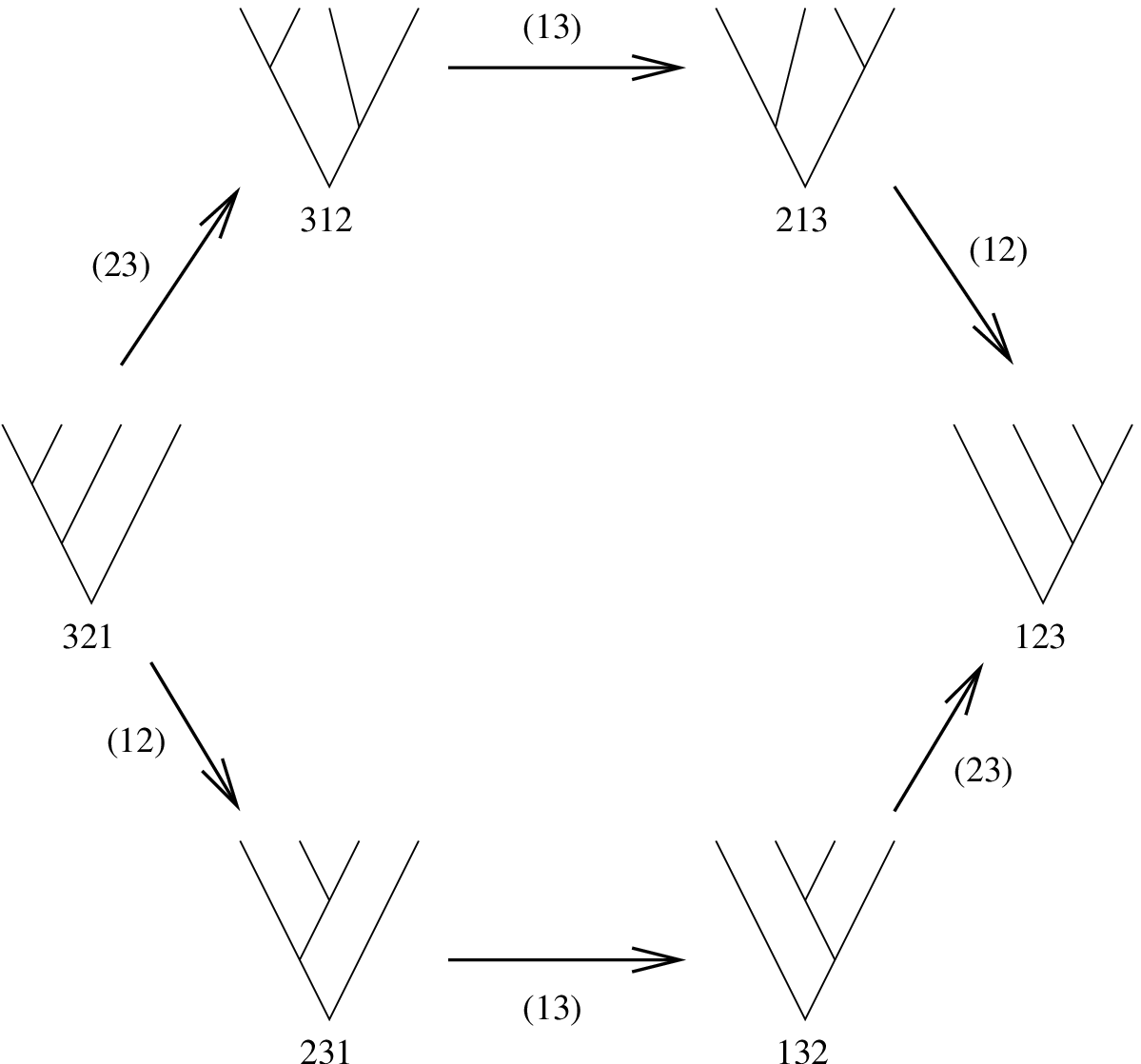}\\ \\
and maps under the functor $\Gamma $ to the following diagram in $\Funct
(\cC^{4},\cC )$.\\ \\
\def\objectstyle{\scriptstyle}
\xymatrix@R=60pt@C=25pt{ & \otimes (\otimes \times \otimes )
  \ar[r]^*+{\qss } &
  \otimes (\otimes \times \otimes ) \ar[dr]^*+{\ass (1\times 1\times
    \otimes )} & \\
  \otimes (\otimes \times 1)(\otimes \times 1\times 1) \ar[ur]^*+{\ass
    (\otimes \times 1\times 1)} \ar[dr]_*+{\otimes (\ass \times 1)} &
  & & \otimes (1\times \otimes )(1\times 1\times \otimes ) \\
  & \otimes (\otimes \times 1)(1\times \otimes \times 1) \ar[r]_*+{\ass
    (1\times \otimes \times 1)} &
  \otimes (1\times \otimes )(1\times \otimes \times 1)
  \ar[ur]_*+{\otimes (1\times \ass )} &
  \\ }\\

Next we introduce some important operations that may be performed on
coupling trees. Let $t$ be a coupling tree of length $n+1$. We define $\{
t\} =[n]$. Every tree defines a partial ordering on $\{ t\} $ where
$i\leq_{t}j$ if level $i$ is connected to level $j$ without passing
through levels less than $i$. In particular $1\leq_{t}i$ for all levels
$i\in \{ t\}$.

We define the cut operations to be maps $\wedge_{i},\vee_{j}:(\Cptr
)_{0}\rightarrow (\Cptr )_{0}$ where $i,j\in \bN $ defined as
follows.  Given a coupling tree $t$, if one cuts the tree at the branch
at level $i$ (if it exists) then one obtains two trees. The upper tree
gives rise to a coupling tree denoted $\vee_{i}t$ while the bottom gives 
rise to a coupling tree denoted $\wedge_{i}t$. If no branch exists at
level $i$ then $\vee_{i}t=0$ is the null tree, and $\wedge_{i}t=t$.
\begin{lem}
The cut operations satisfy the following properties:
\begin{enumerate}
\item{$\wedge_{1}t=1$ and $\wedge_{i}t=t$ for
    all $i\geq |t|$ and all coupling trees $t$.}
\item{$\vee_{1}t=t$, $|\vee_{|t|-1}t|=2$, $\vee_{i}t=0$ for all
    $i\geq |t|$ and all coupling trees $t$.}
\item{$\wedge_{i}=\wedge_{i}\wedge_{j}$ whenever $i\leq_{t}j$. In
    particular $\wedge_{i}$ is idempotent.}
\item{Given a coupling tree $t$ then $i\leq_{t}j$ implies
    $\vee_{k}\vee_{i}t=\vee_{j}t$ for some $k\leq |i-j|$.}
\end{enumerate}
\end{lem}
These properties are easily demonstrated.

The reattachment operation $p(n):t\rightarrow t^{\prime }$ at the $n$th
level is an arrow satisfying $n\leq_{t}n+1$ and $n\leq_{t^{\prime }}n+1$
that only interchanges levels $n$ and $n+1$. It is called a reattachment
to the left (resp. right) if $t^{-1}n<t^{-1}(n+1)$ (resp.
$t^{-1}n>t^{-1}(n+1)$). Hence $\vee_{n+1}t=\vee_{n+1}t^{\prime }$ and
$\wedge_{n}t=\wedge_{n^{\prime }}t$. For a left reattachment
$R\wedge_{n}t=\wedge_{n+1}t$, and for a right reattachment
$L\wedge_{n}t=\wedge_{n+1}t$.
\begin{defn}
  An arrow of $\Cptr $ is called primitive if it is an identity or
  corresponds to a single reattachment operation.
\end{defn}

{\it Proof of Theorem \ref{cohthm1}:} First let $\rho (m):t\rightarrow
t^{\prime }$ be a reattachment arrow at the $m$th level. If $\rho (m)$
reattaches to the right we define
\begin{eqnarray}
\Gamma \rho (m) & = & \Gamma \wedge_{m}t \left( 1^{p}\times \ass
  (\Gamma LL\vee_{m}t\times \Gamma RL\vee_{m}t\times \Gamma
  R\vee_{m}t)\right. \nonumber \\
& & \ph{spacespacespace}\left. \times 1^{q}\right)
\end{eqnarray}
where $p+1=\mini \{ t^{-1}n:m\leq_{t}n\} $ and
$q=|t|-p-|\vee_{m}t|$. Otherwise it reattaches to the left and we define
\begin{eqnarray}
\Gamma \rho (m) & = & \Gamma \wedge_{m}t \left( 1^{p}\times \ass^{-1}
  (\Gamma L\vee_{m}t\times \Gamma LR\vee_{m}t\times \Gamma
  RR\vee_{m}t)\right. \nonumber \\
& & \ph{spacespacespace}\left. \times 1^{q}\right)
\end{eqnarray}
For any arrow $f$ of $\Cptr $ there is a (directed) sequence of
primitive reattachment arrows with $f=\rho_{n}\cdots \rho_{1}$. Then we
define $\Gamma (f)=(\Gamma \rho_{n})\cdots (\Gamma \rho_{1})$. It only
remains to show that $\Gamma $ is well--defined and a functor.
Equivalently we show that any commutative polygonal diagram of
primitives in $\Cptr $ maps to a commutative diagram in $\Funct(
\cC^{\bW },\cC )$. Proof is by induction on coupling tree length $n$. It
is easy to show for $n=1,2,3,4$.  Consider a polygonal diagram in $\Cptr
$ with vertices $t_{0},...,t_{N-1}$ of length $n+1>4$, where we take
$t_{0}=t_{N}$, and all arrows are primitive. We show that the
commutativity of this diagram is equivalent to the commutativity of a
diagram in which the highest level $q=n-1$ is maintained in a fixed
region for all of its vertices.  Such a diagram is commutative by the
induction hypothesis. Let $r=\maxi \{t_{k}^{-1}q:k\in [N]\} $ be the
right--most region containing the level $q$ in some vertex. Suppose
$t_{k}\rightarrow \cdots \rightarrow t_{l}$ is a section of the diagram
where $t^{-1}_{k-1}q,t^{-1}_{l+1}q<r$ and $t^{-1}_{i}q=r$ whenever
$k\leq i\leq l$. We replace this section with an alternative section
$t_{k-1}\rightarrow s_{1}\rightarrow \cdots \rightarrow s_{d}\rightarrow
t_{l+1}$ such that the enclosed region commutes. Iterating this
procedure until level $q$ is in a fixed region will complete the proof.
The arrows $t_{k-1}\rightarrow t_{k}$ and $t_{l}\rightarrow t_{l+1}$ are
primitives at the $q$th level. There exists a primitive sequence
$t_{k}\rightarrow u_{1}\rightarrow \cdots \rightarrow u_{d}\rightarrow
t_{l}$ keeping the levels $q$ and $q-1$ fixed. Moreover the enclosed
diagram under $\Gamma $ commutes by the induction hypothesis since level
$q$ is kept fixed. Next construct the same sequence of operations
starting with $t_{k-1}$ giving a sequence $t_{k-1}\rightarrow
s_{1}\rightarrow \cdots \rightarrow s_{d}\rightarrow t_{k+1}$. This
encloses a diagram with the previous sequence that commutes under
$\Gamma $ because it is a ladder of natural squares. This is the desired
replacement sequence completing the proof.

For a premonoidal structure the primitive reattachment operations are
restricted to adjacent levels. In the pseudo--monoidal situation the
adjacent restriction is lifted. Thus a primitive arrow $\rho
(n):t\rightarrow t^{\prime }$ for (pseudo--monoidal) reattachment at the
$n$th level satisfies $q\equiv \mini \{ m:n<_{t}m\} =\mini \{
m:n<_{t^{\prime }}m\} $ and only interchanges levels $n$ and $q$.  We
say an arrow $\pi :s\rightarrow t$ is split about the $m$ level if $\{
l:m\leq_{s}l\} =\{ l:m\leq_{t}l\} $. Thus we can write $\pi =\sigma \tau
=\tau \sigma $ where $\sigma $ permutes only the levels in $\{
l:m<_{s}l\} $ and $\tau $ only the levels in $\{ l:m\nleq_{m}l\} $. We
have the following coherence result from Joyce \cite{wj}.
\begin{thm} \label{cohthm2}
If $(\otimes ,\ass )$ is a pseudo--monoidal structure then $\Gamma :\Cptr
\rightarrow \Funct (\cC^{\bW },\cC )$ of theorem \ref{cohthm1} satisfies
\begin{eqnarray}
\Gamma (\sigma \tau ) & = & \otimes (\Gamma (\sigma ) \times \Gamma
(\tau )) \label{propty2}\\
\Gamma (1) & = & 1
\end{eqnarray}
for any arrow of $\Cptr $ split about level $2$ into $\sigma $ and $\tau $.
\end{thm}
{\it Proof } The proof mirrors that of theorem \ref{cohthm1} except we
have a different procedure for calculating the alternative sequence
$t_{k-1}\rightarrow s_{1}\rightarrow \cdots \rightarrow s_{d}\rightarrow
t_{l+1}$ in the induction step. Suppose $t_{k-1}\rightarrow t_{k}$ is a
reattachment at the $m$th level. We replace this arrow with the
following sequence. Suppose $t^{-1}(m+1)<t^{-1}m$ (the other case is
shown similarly and left to the reader) and consider the following
diagram.\\ \\
\hspace*{2cm}\includegraphics[width=150pt]{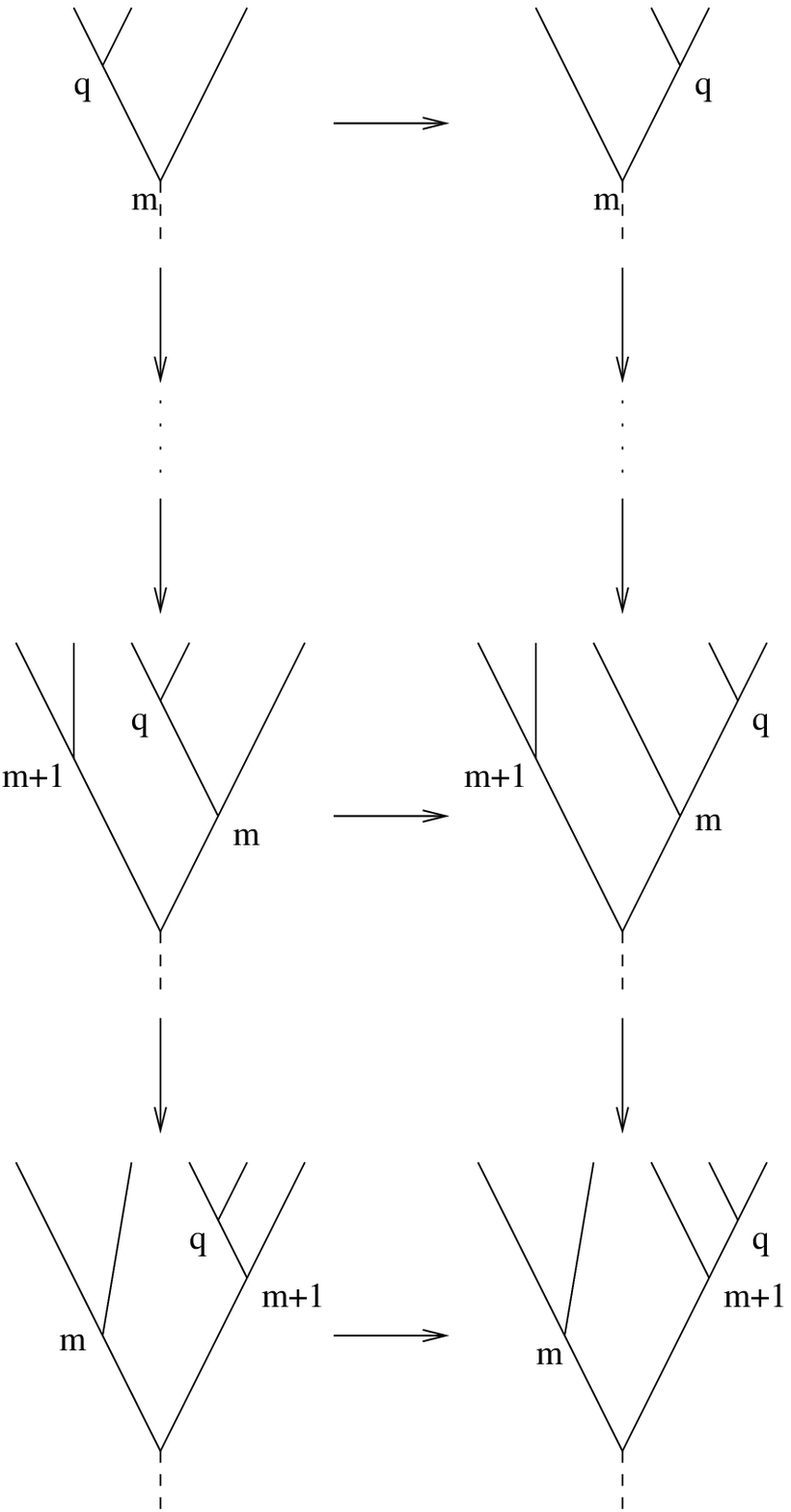}\\ \\
The top arrow is $t_{k-1}\rightarrow t_{k}$. The sides of the top
diagram are parallel operations keeping the subtree with root $m$ fixed.
The bottom diagram is a $\qss $--square. Hence substituting for
$t_{k-1}\rightarrow t_{k}$ using the outside sequence we re--identify
the maximal sequence with $t_{k-1}\rightarrow t_{k}$ reattaching at
level $m+1$. Inductively we are led to a maximal sequence with
$t_{k-1}\rightarrow t_{k}$ reattaching at level $q-1$. Similarly we are
led to $t_{l}\rightarrow t_{l+1}$ reattaching about level $q-1$. We can
construct an alternative sequence from $t_{k}$ to $t_{l}$ keeping $q$
and $q-1$ fixed. The diagram enclosed commutes by the induction
hypothesis. If we apply the same sequence of operations from $t_{k-1}$
to $t_{l+1}$ we enclose a ladder of natural squares. Moreover, this is
the desired replacement sequence completing the proof.

\section{Coupling Trees and Monoidal Structures without Unit}        %

Underlying every coupling tree is a planar rooted binary tree, or
bracketing, given by simply forgetting the levels. This allows us to
define an equivalence relation where $s\sim t$ if and only if both trees
give the same bracketing under forgetting levels. Thus we have a full
forgetful functor $U:\Cptr \rightarrow \Cptr /\sim $ onto the quotient
category $\Cptr /\sim $.
\begin{prop} \label{prop1}
The forgetful functor $U:\Cptr \rightarrow \Cptr /\sim $ given by forgetting 
levels has a right adjoint.
\end{prop}
Let $[t]=\{ s:s\sim t\} $ and let $\phi :(\Cptr )_{0}/ \sim\ \rightarrow
(\Cptr )_{0}$ be a choice function choosing an isotypical member for
each equivalence class. Thus $[\phi [t]]=[t]$. Define the faithful
functor $M:\Cptr /\sim\ \rightarrow \Cptr $ given by $\phi $ on objects
assigning the unique arrow between any two objects of the same length.
Clearly $UM=1$ so the counit is strict and the unit $\eta :MU\rightarrow
1$ is given by assigning $\eta_{t}$ to be the unique arrow $t\rightarrow
\phi [t]$.

When making a choice of isotypical objects without invoking the axiom of
choice one needs a criterion. For example, define the following order on
coupling trees. Given $s$ and $t$ we define $s<t$ if $|s|<|t|$. If
$|s|=|t|$ then $s<t$ if there is some $k$ such that $sj=tj$ for all
$j<k$ and $sk<tk$. Thus $\phi [t]$ can be chosen to be the maximal
(alternatively minimal) member of $[t]$.

There is no canonical choice of tensor product on $\Cptr $, however,
there is on $\Cptr /\sim $. Given two bracketings $b_{1}$ and $b_{2}$
there is a unique bracketing $b_{3}$ such that $ULMb_{3}=b_{1}$ and
$URMb_{3}=b_{2}$. We then define $b_{1}\otimes b_{2}=b_{3}$ which extends
to a unique bifunctor $\otimes :\Cptr /\sim \times\ \Cptr /\sim\ 
\rightarrow \Cptr /\sim $. This defines a unique monoidal structure on
$\Cptr /\sim $. The adjunction $U\dashv M$ of proposition \ref{prop1}
can be used to lift this bifunctor to $\Cptr $ by defining $\otimes_{M}
\equiv M\otimes (U\times U):\Cptr \times \Cptr \rightarrow \Cptr $. Thus
each $\otimes_{M}$ admits a unique monoidal structure $(\otimes_{M}
,\ass^{M}\equiv M\ass (U\times U\times U))$ on $\Cptr $.
Unlike for $\otimes $ on $\Cptr /\sim $, no single bifunctor $\otimes_{M}$
generates all objects of $\Cptr $ from a single generator.

\begin{prop} \label{prop2}
  Let $\cC $ be a category with premonoidal structure $(\otimes ,\ass
  )$. The functor $\Gamma :\Cptr \rightarrow \Funct (\cC^{\bW },\cC )$
  of theorem \ref{cohthm1} is such that, for all coupling trees $s,t,u$
  we can find $p:(s\otimes_{M}t)\otimes_{M}u)\rightarrow
  M[(s\otimes_{M}t)\otimes_{M} u]$ and $q:s\otimes_{M}(t\otimes_{M}u)
  \rightarrow M[s\otimes_{M}(t\otimes_{M}u)]$ such that the following
  square commutes.\\ \\
  \xymatrix@R=40pt@C=60pt{ \Gamma \left( (s\otimes_{M}t)\otimes_{M}u\right)
    \ar[r]^{\Gamma \ass^{M}_{s,t,u}}
    \ar[d]_{\Gamma p} & \Gamma \left( s\otimes_{M}(t\otimes_{M}u)\right) \\
    \otimes (\otimes \times 1)(\Gamma s\times \Gamma t\times \Gamma u)
    \ar[r]_{\ass (\Gamma s\times \Gamma t\times \Gamma u)} & \otimes
    (1\times \otimes )(\Gamma s\times \Gamma t\times \Gamma u)
    \ar[u]_{\Gamma q}\\ }
\end{prop}

The notion of a premonoidal functor is that of a monoidal functor
without the properties pertaining to the unit.
\begin{defn}
  Given categories $\cC $ and $\cC^{\prime }$ with premonoidal
  structures $(\otimes ,\ass )$ and $(\otimes^{\prime },\ass^{\prime })$
  respectively, a premonoidal functor is a pair $(F,\phi )$ where $F:\cC
  \rightarrow \cC^{\prime }$ is a functor and $\phi :\otimes^{\prime
    }(F\times F)
  \rightarrow F\otimes $ is a natural transformation satisfying\\ \\
  \xymatrix@R=40pt@C=60pt{ \otimes^{\prime } (\otimes^{\prime }\times
    1)(F\times F\times F)\ar[r]^{\ass^{\prime }(F\times F\times F)}
    \ar[d]_{\otimes^{\prime }(\phi \times F)} & \otimes^{\prime
      }(1\times \otimes^{\prime
      })(F\times F\times F)\ar[d]^{\otimes^{\prime }(F\times \phi )}\\
    \otimes^{\prime }(F\otimes \times F) \ar[d]_{\phi (\otimes \times
      1)} &
    \otimes^{\prime }(F\times F\otimes ) \ar[d]^{\phi (1\times \otimes )}\\
    F\otimes (\otimes \times 1) \ar[r]_{F\ass } & F\otimes (1\times \otimes )\\
    }\\ \\
  A premonoidal functor $(F,\phi )$ is called strong (resp. strict) if
  $\phi $ is an isomorphism (resp. identity).
\end{defn}
We have a (weakened) restatement of Mac Lane's coherence theorem
\cite{sm2} for monoidal categories without a unit.
\begin{cor}
  If $\cC $ is a category with a monoidal structure without unit
  $(\otimes ,\ass )$ then $\Gamma M:\Cptr /\sim\ \rightarrow \Funct
  (\cC^{\bW },\cC )$ is a strong premonoidal functor.
\end{cor}
This follows from noting that $M:\Cptr /\sim\ \rightarrow \Cptr $ is a
strong monoidal functor and applying proposition \ref{prop2} where
$\Gamma p=1$ and $\Gamma q=1$ for a monoidal structure.

\section{Premonoidal Structures with Unit}                           %

In this section we add a unit without introducing triangle constraints.
\begin{defn}
  A premonoidal (resp. pseudo--monoidal) structure with unit for a
  category $\cC $ is a pentuple $(\otimes ,\ass ,\lid ,\rid ,e)$ where
  $(\otimes ,\ass )$ is a premonoidal (resp. pseudo--monoidal) structure
  for $\cC $, $e$ is an object of $\cC $ called the unit object, and
  $\lid :\otimes (I\times 1)\rightarrow \pi_{2}$ and $\rid :\otimes
  (1\times I)\rightarrow \pi_{2}$ are natural isomorphisms called
  respectively the recouplings for left unit and right unit.
\end{defn}
The functor $I:\cC \rightarrow \cC $ is defined by $If=1_{e}$ for all
arrows $f$ and may be called the unit functor. The functors
$\pi_{1},\pi_{2}:\cC^{2}\rightarrow \cC $ are given by the universal
projections of the Cartesian product. That is $\pi_{k}(f_{1},f_{2})=f_{k}$
for all arrows $(f_{1},f_{2})$ of $\cC^{2}$ and $k=1,2$. Thus for any arrow
$(f,g):(a,b)\rightarrow (c,d)$ the recouplings for left unit and right
unit satisfy the natural squares\\ \\
\xymatrix@R=40pt@C=40pt{
e\otimes b \ar[r]^{\lid_{e,b}} \ar[d]_{1_{e}\otimes g} & b
    \ar[d]^{g} \\ e\otimes d \ar[r]_{\lid_{e,d}} & d \\
}\hspace*{2cm}
\xymatrix@R=40pt@C=40pt{
a\otimes e \ar[r]^{\rid_{a,e}} \ar[d]_{f\otimes 1_{e}} & a \ar[d]^{f} \\ 
c\otimes e \ar[r]_{\lid_{c,e}} & d \\
}\\ \\
The labeling of the left and right unit recouplings are redundant so often
we identify $\lid_{b}\equiv \lid_{e,b}$ and $\rid_{a}\equiv
\rid_{a,e}$. In fact this identification gives the usual form of the
left and right unit recouplings as $\lid :e\otimes \ubar \rightarrow 1$
and $\rid :\ubar \otimes e\rightarrow 1$.

The triangle diagrams do not hold so we measure their non--commutativity
by defining {\it ghost} natural automorphisms.
\begin{defn}
We define the ghosts for associativity to be the natural automorphisms
$\gss (12) ,\gss (23) ,\gss (13):\otimes \rightarrow \otimes $ defined by
\begin{eqnarray}
\gss (23) (\pi_{2}\times 1) & = & \left(\lid (1\times \otimes )\right)
\left( \ass (I\times 1\times 1)\right) \left( \otimes (\lid^{-1}\times
  1)\right) \\
\gss (13) (\pi_{1}\times 1) & = & \left( \otimes (1\times \lid )\right)
\left( \ass (1\times I\times 1)\right) \left( \otimes (\rid^{-1}\times
  1)\right) \\
\gss (12) (1\times \pi_{1}) & = & \left( \otimes (1\times \rid )\right)
\left( \ass (1\times 1\times I)\right) \left( \rid^{-1} (\otimes \times
  1)\right)
\end{eqnarray}
\end{defn}
The ghost associativity natural automorphisms satisfy (and are defined
by) the following ghostly triangle diagrams.\\ \\
\xymatrix@R=40pt@C=40pt{
(e\otimes b)\otimes c \ar[r]^{\ass_{e,b,c}} \ar[d]_{\lid_{b}\otimes
  1_{c}} & e\otimes (b\otimes c) \ar[d]^{\lid_{b\otimes c}} \\
b\otimes c \ar[r]_{\gss (23)_{b,c}} & b\otimes c \\ }\hspace*{2cm}
\xymatrix@R=40pt@C=40pt{
(a\otimes e)\otimes c \ar[r]^{\ass_{a,e,c}} \ar[d]_{\rid_{a}\otimes
  1_{c}} & a\otimes (e\otimes c) \ar[d]^{1_{a}\otimes \lid_{c}} \\
a\otimes c \ar[r]_{\gss (13)_{a,c}} & a\otimes c \\ }\\ \\
\hspace*{4cm}\xymatrix@R=40pt@C=40pt{
(a\otimes b)\otimes e \ar[r]^{\ass_{a,b,e}} \ar[d]_{\rid_{a\otimes b}} &
a\otimes (b\otimes e) \ar[d]^{1_{a}\otimes \rid_{b}} \\
a\otimes b \ar[r]_{\gss (12)_{a,b}} & a\otimes b \\ }\\ \\
for all objects $a,b,c$ of $\cC $.

Note that the associative structure with unit is monoidal whenever the
deformativity and ghost natural automorphisms are identities. The
ghostly triangle diagram for $\gss (13)$ becomes the triangle
constraint, and together with the Pentagonal constraint imply the other
triangle constraints \cite{gk}. In the monoidal situation, Mac Lane
\cite{sm2} has proved a well--known coherence result.

Similarly we can define ghosts for deformativity $\gss (234),\gss
(134):\otimes (1\times \otimes )\rightarrow \otimes (1\times \otimes )$
and $\gss (124),\gss (123):\otimes (\otimes \times 1)\rightarrow \otimes 
(\otimes \times 1)$ according to
\begin{eqnarray}
\gss (234)(\pi_{2}\times 1\times 1) & = & \left( \otimes (\lid \times
\otimes )\right) \left( \qss (I\times 1\times 1\times 1)\right) \left(
\otimes (\lid^{-1}\times \otimes )\right) \nonumber \\
\gss (134)(\pi_{1}\times 1\times 1) & = & \left( \otimes (\rid \times
\otimes )\right) \left( \qss (1\times I\times 1\times 1)\right) \left(
\otimes (\rid^{-1}\times \otimes )\right) \nonumber \\
\gss (124)(1\times 1\times \pi_{2}) & = & \left( \otimes (\otimes
\times \lid )\right) \left( \qss (1\times 1\times I\times 1)\right) \left(
\otimes (\otimes \times \lid^{-1})\right) \nonumber \\
\gss (123)(1\times 1\times \pi_{1}) & = & \left( \otimes (\otimes
\times \rid )\right) \left( \qss (1\times 1\times 1\times I)\right) \left(
\otimes (\otimes \times \rid^{-1})\right) \nonumber
\end{eqnarray}
We leave it to the reader to write down diagrams.

\section{Coherence of Premonoidal Structures with Unit}              %

We extend the groupoid of coupling trees by attaching two types of
nodules on leaves, called unit and ghost nodules. Given a finite set $U$ we
construct a groupoid $\cN (U)$ called the nodule groupoid over $U$. The
objects are pairs $(u,v)$ where $u,v\subset U$, $u\cap v=\emptyset $ and
$v\neq U$. We say the object $(u,v)$ represents $|u|$ unit nodules and $|v|$
ghost nodules. There is at most one arrow between two objects given by
the condition: $(u,v)\rightarrow (u^{\prime },v^{\prime })$ is an arrow
if and only if $u\cup v=u^{\prime }\cup v^{\prime }$. In other words the 
arrows interchange nodule type. We see that
\begin{eqnarray}
\cN (U) & = & \coprod_{k=0}^{|U|} \cN_{k}(U)
\end{eqnarray}
where $\cN_{k}(U)$ is the full subgroupoid whose objects are given by
$(u,v)$ such that $|u|+|v|=k$. The groupoid of coupling trees with
nodules is given by
\begin{eqnarray}
\NCptr & = & \coprod_{n\in \bW }\Cptr_{n}\times \cN ([n])\\
& = & \coprod_{n\in \bW }\coprod_{k=0}^{n}\Cptr_{n}\times \cN_{k}([n])
\end{eqnarray}
A noduled coupling tree $(t,u,v)$ has $|t|$ leaves, with unit
nodules (open circles) in positions $i$ for all $i\in u$, and ghost
nodules (closed circles) in positions $j$ for all $j\in v$. For
example $(514632,\{ 3\} ,\{ 5,6\} )$ is given by\\ \\
\hspace*{2cm}
\includegraphics[width=100pt]{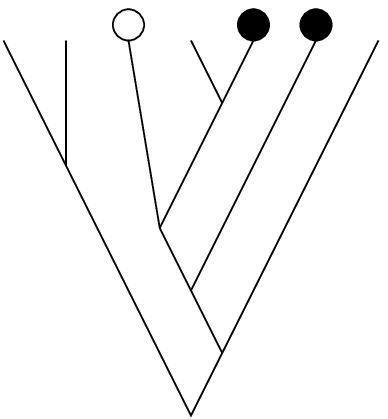}\\ \\
The left and right coupling tree operations are extended to noduled
coupling trees by
\begin{eqnarray}
L(t,u,v) & = & (Lt,u\cap [m],v\cap [m])\\
R(t,u,v) & = & (Rt,u\setminus [m]-m,v\setminus [m]-m)
\end{eqnarray}
provided $v\cap [m]\neq [m]$ and $(v\setminus [m]-m)\cap [|Rt|]\neq
[|Rt|]$ where $m=|Lt|$ and we have defined $U+k\equiv \{ i+k:i\in U\} $
for $U\subset \bZ $ and $k\in \bZ $.

We make a few convenient definitions. Define $1^{k}=1\times 1 \times
\cdots \times 1:\cC^{k}\rightarrow \cC^{k}$ and
$\pi^{k}_{i}:\cC^{k}\rightarrow \cC $, the latter taking
$(f_{1},...,f_{k})\mapsto f_{i}$. Also given a set $U\subset \bN $ and a
coupling tree $t$ we define $C_{U}t$ to be the unique coupling tree
obtained by contracting out those leaves whose positions are in the set
$U$.

We define an equivalence relation on $\NCptr $ as follows. We write
$(s,u,v)\sim (t,w,x)$ if and only if $s\sim t$, $C_{v}s=C_{x}t$, $u=w$
and $v=x$. This equivalence extends uniquely to arrows. This defines a
forgetful functor $U:\NCptr \rightarrow \NCptr /\sim $ determined by
mapping each arrow to its equivalence class. As in proposition
\ref{prop1}, $U$ has right adjoint sections $M$ given by choosing a
representative member of each equivalence class of objects. Each
equivalence class may be thought of as the coupling tree without levels
and the edges attached to ghost nodules omitted as in the following
example.\\ \\
\hspace*{2cm}
\includegraphics[width=200pt]{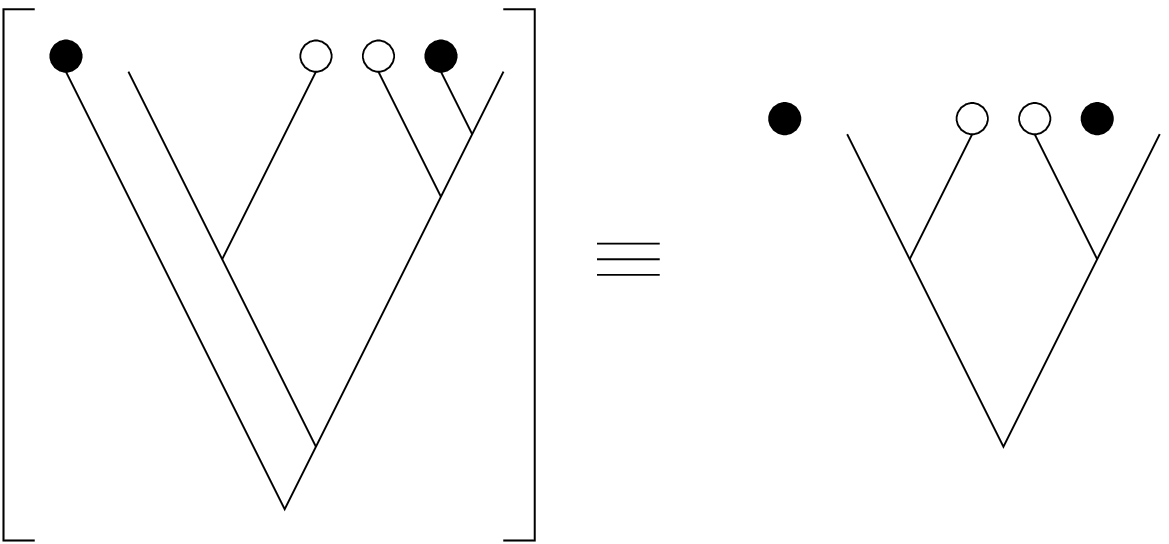}\\ \\
The category $\NCptr /\sim $ is a monoidal category with bifunctor
$\otimes $ given by concatenation and joining the roots to the leaves of
the tree $1$.  This can be lifted to a bifunctor on $\NCptr $ as
$\otimes_{M}=M\otimes (U\times U)$.

\begin{defn}
Given categories $\cC $ and $\cC^{\prime }$ with monoidal structures
$(\otimes ,\ass ,\lid ,\rid ,e)$ and $(\otimes^{\prime },\ass^{\prime }, 
\lid^{\prime },\rid^{\prime },e^{\prime })$
respectively, a monoidal functor is a triplet $(F,\phi ,\psi )$ where
$(F,\phi )$ is a premonoidal functor and $\psi :FI\rightarrow I^{\prime
  }F$ a natural transformation satisfying\\ \\
\xymatrix@R=40pt@C=60pt{
\otimes^{\prime }(I^{\prime }F\times F) \ar[r]^{\lid^{\prime }(F\times
  F)} & F\pi_{2} \\  \otimes^{\prime }(FI\times F)
\ar[u]^{\otimes^{\prime }(\psi \times F)} \ar[r]_{\phi (I\times 1)} &
  F\otimes (I\times 1) \ar[u]_{F\lid }
}\xymatrix@R=40pt@C=60pt{
\otimes^{\prime }(F\times I^{\prime }F) \ar[r]^{\rid^{\prime }(F\times
  F)} & F\pi_{1} \\  \otimes^{\prime }(F\times FI)
\ar[u]^{\otimes^{\prime }(F\times \psi )} \ar[r]_{\phi (1\times I)} &
  F\otimes (1\times I) \ar[u]_{F\rid }
}\\ \\
A monoidal functor $(F,\phi ,\psi )$ is called strong (resp. strict) if
$\phi $ and $\psi $ are isomorphisms (resp. identities).
\end{defn}

We extend the functor of theorem \ref{cohthm1} to $\NCptr $ giving the
following coherence result.
\begin{thm} \label{cohthm3}
  Given a category $\cC $ and a premonoidal structure with unit
  $(\otimes ,\ass ,\lid ,\rid ,e)$ there is an extension of $\Gamma
  :\Cptr \rightarrow \cC $ to $\NCptr $ such that the arrows $(1,\{ 1\}
  ,\emptyset )\rightarrow (1,\emptyset ,\{ 1\} )$ and\\ $(1,\{ 2\}
  ,\emptyset )\rightarrow (1,\emptyset ,\{ 2\} )$ map under $\Gamma $ to
  $\lid $ and $\rid $ respectively. Furthermore, if $(\otimes ,\ass )$
  is pseudo--monoidal and all the ghosts vanish then $\Gamma M:\NCptr
  /\sim\ \rightarrow \cC $ is a monoidal functor for all $\otimes_{M}$.
\end{thm}

{\it Proof} $\Gamma $ is characterised inductively on objects as
follows. We define
\begin{eqnarray}
\Gamma (t,\emptyset ,\emptyset ) & = & \Gamma t\\
\Gamma (\emptyset ,\{ 1\} ,\emptyset ) & = & e
\end{eqnarray}
Whenever $\Gamma L(t,u,v)$ and $\Gamma R(t,u,v)$ are defined then
\begin{eqnarray}
\Gamma (t,u,v) & = & \otimes (\Gamma L(t,u,v)\times \Gamma R(t,u,v)) ~,
\end{eqnarray}
If $|\Gamma Lt|=2$ then
\begin{eqnarray}
\Gamma (t,u,v) & = & \Gamma R(t,u,v)(\pi_{i}\times 1^{|Rt|})
\end{eqnarray}
whenever $\{ i\} =v$, and if $|Rt|=2$ and then
\begin{eqnarray}
\Gamma (t,u,v) & = & \Gamma L(t,u,v)( 1^{|Lt|}\times \pi_{i})
\end{eqnarray}
whenever $\{ i\} =v$. We take the primitive arrows $s\rightarrow t$ to
be either reattachment arrows at the $n$th level where each of
$L\vee_{n+1}s$, $R\vee_{n+1}s$, $L\vee_{n}s$ and $R\vee_{n}s$ contain a
ghost nodule free leaf, or to be nodule change arrows where a single
nodule type is changed. The image of a reattachment arrow under $\Gamma
$ is given by theorem \ref{cohthm1}. Let $\rho (m):t\rightarrow t^{\prime
  }$ be a nodule change arrow converting a unit nodule into a ghost
nodule at the $m$th level (a nodule in position $t^{-1}m$ is
changed). Let $p+1=\mini \{ t^{-1}n:m\leq_{t}n\} $ and
$q=p+|\vee_{m}t|$. If $t^{-1}m=p+1$ we define
\begin{eqnarray}
\Gamma \rho (m) & = & \left( \Gamma \wedge_{m}t \right) \left(
  1^{p}\times \lid (I\times \Gamma R\vee_{m}t)\times 1^{q}\right) \nonumber
\end{eqnarray}
Otherwise $t^{-1}m=q-1$ and we define
\begin{eqnarray}
\Gamma \rho (m) & = & \left( \Gamma \wedge_{m}t\right) \left(
  1^{p}\times \rid (\Gamma L\vee_{m}t\times I)\times 1^{q}\right) \nonumber
\end{eqnarray}
For any arrow $f$ of $\NCptr $ there is a (directed) sequence of
primitive reattachment arrows with $f=\rho_{n}\cdots \rho_{1}$. We then
define $\Gamma (f)=(\Gamma \rho_{n})\cdots (\Gamma \rho_{1})$. It only
remains to show that $\Gamma $ is well--defined and a functor.
Equivalently we show that any commutative diagram of primitives in
$\NCptr $ maps to a commutative diagram in $\Funct( \cC^{\bW },\cC )$.
Specifically we show how to remove ghost nodules in the $i$th position.
Then all ghost nodules can be removed and the result follows from
theorem \ref{cohthm1}. It is not hard to see that the primitive arrows
for changing ghost nodules into unit nodules commute with all other
primitive arrows. Hence the sections of the diagram with ghost nodules
in the $i$th position can be replaced by an alternative sequence without
ghost nodules enclosing a ladder of natural squares.

\section{Braidings}                                                  %

We introduce a braid structure on a category requiring only that it
possess a premonoidal structure.
\begin{defn}
A prebraid structure $(\otimes ,\ass ,\com )$ for a category $\cC $ is
a premonoidal structure $(\ass ,\otimes )$ and a natural isomorphism
$\com :\otimes \rightarrow \otimes \tau_{(12)}$ where $\tau_{(12)}$ is
the switch match. This structure is called braid premonoidal if the
following three diagrams\\ \\
\xymatrix@R=40pt@C=40pt{ (a\otimes b)\otimes c \ar[r]^{\com_{a\otimes
        b,c}} \ar[d]_{\ass_{a,b,c}} & c\otimes (a\otimes b) \\
    a\otimes (b\otimes c) \ar[d]_{1_{a}\otimes \com_{b,c}} & (c\otimes
    a)\otimes b \ar[u]_{\ass_{c,a,b}} \\
    a\otimes (c\otimes b) \ar[r]_{\ass^{-1}_{a,c,b}} & (a\otimes
    c)\otimes b \ar[u]_{\com_{a,c}\otimes 1_{b}} \\ }
\hspace*{1cm} \xymatrix@R=40pt@C=40pt{ a\otimes
    (b\otimes c) \ar[r]^{\com_{a,b\otimes c}} \ar[d]_{\ass^{-1}_{a,b,c}}
    & (b\otimes c)\otimes a \\
    (a\otimes b)\otimes c \ar[d]_{\com_{a,b}\otimes 1_{c}} & b\otimes
    (c\otimes a) \ar[u]_{\ass^{-1}_{b,c,a}} \\
    (b\otimes a)\otimes c \ar[r]_{\ass_{b,a,c}} & b\otimes (a\otimes c)
    \ar[u]_{1_{b}\otimes \com_{a,c}}
    \\ }\\ \\ \\
\xymatrix@C=40pt@R=40pt{ (a\otimes b)\otimes (c\otimes d)
    \ar[r]^{\qss_{a,b,c,d}} \ar[d]_{\com_{a\otimes b,c\otimes d}} &
    (a\otimes b)\otimes
    (c\otimes d) \ar[d]^{\com_{a\otimes b,c\otimes d}} \\
    (c\otimes d)\otimes (a\otimes b) & (c\otimes
    d)\otimes (a\otimes b) \ar[l]^{\qss_{c,d,a,b}} \\ }\\ \\ \\
  commute for all objects $a,b,c,d$ of $\cC $. If in addition $(\otimes
  ,\ass )$ is pseudo--monoidal and the following square diagram\\ \\
\xymatrix@C=40pt@R=40pt{ (a\otimes b)\otimes (c\otimes d)
    \ar[r]^{\qss_{a,b,c,d}} \ar[d]_{\com_{a,b}\otimes 1_{c\otimes d}} &
    (a\otimes b)\otimes
    (c\otimes d) \ar[d]^{\com_{a,b}\otimes 1_{c\otimes d}} \\
    (b\otimes a)\otimes (c\otimes d) \ar[r]_{\qss_{b,a,c,d}} & (b\otimes
    a)\otimes (c\otimes d) \\ }\\ \\ \\
  commutes for all objects $a,b,c,d$ of $\cC $ then the structure is
  called braid pseudo--monoidal. Finally, whenever $\com^{-1}=\com
  \tau_{(12)}$ the braid is called a symmetry.
\end{defn}
More generally the switch map extends to an action $\tau
:S_{n}\rightarrow \End (\cC^{n})$ where $\pi \mapsto \tau_{\pi }$ is
given by $\tau_{\pi }(c_{1},...,c_{n})=(c_{\pi 1},...,c_{\pi n})$.

\begin{defn}
Given categories $\cC $ and $\cC^{\prime }$ with prebraid structures
$(\otimes ,\ass ,\com )$ and $(\otimes^{\prime },\ass^{\prime }, 
\com^{\prime })$ respectively, a braid premonoidal functor is a
premonoidal functor $(F,\phi )$ satisfying\\ \\
\xymatrix@R=40pt@C=60pt{
\otimes^{\prime }(F\times F) \ar[r]^{\com^{\prime }(F\times F)}
\ar[d]_{\phi } & \otimes^{\prime }(F\times F)\tau_{(12)} \ar[d]^{\phi
  \tau_{(12)}} \\
F\otimes \ar[r]_{F\com } & F\otimes \tau_{(12)}
}\\ \\
A braid premonoidal functor $(F,\phi )$ is called strong (resp. strict)
if $\phi $ is an isomorphism (resp. identity).
\end{defn}

Coherence will be described with respect to the Artin braid groups
$B_{n}$ where $n\in \bW $. These are groupoids on one object. The group 
$B_{n}$ is generated by $\tau_{1},...,\tau_{n-1}$ satisfying the
conditions
\begin{eqnarray}
\tau_{i}\tau_{j} & = & \tau_{j}\tau_{i} \label{brcond1}\\
\tau_{i+1}\tau_{i}\tau_{i+1} & = & \tau_{i}\tau_{i+1}\tau_{i} \label{brcond2}
\end{eqnarray}
for all $i,j=1,2,...,n-1$ satisfying $|i-j|>1$. Let $V:B_{n}\rightarrow
S_{n}$ be the forgetful functor forgetting the order in which braids
cross. This functor is completely determined on generators by
$V(\tau_{i})=(i\, i+1)$. We define the exploded $n$th braid groupoid
$\xB_{n}$ to be given by the formal collection of arrows $\tau :\pi
\rightarrow V(\tau )\pi $ where $\pi \in S_{n}$ and $\tau \in B_{n}$.
Composition is inherited from $B_{n}$ whenever the source and target
match. The objects are given by $S_{n}$ and the hom--sets by
$\xB_{n}(\pi ,\sigma )=\{ \tau \in B_{n}:V(\tau )\pi =\sigma \} $.

We define the braid groupoid of coupling trees by
\begin{eqnarray}
\BCptr & = & \coprod_{n\in \bW }\Cptr_{n}\times \xB_{n}
\end{eqnarray}
We can now state the main coherence result.
\begin{thm} \label{cohthm4}
  Given a braided premonoidal structure for $\cC $ there is an
  extension of theorem \ref{cohthm1} to $\Gamma :\BCptr \rightarrow
  \Funct (\cC^{\bW },\cC )$ where $\Gamma \tau_{1}=\com $ on $\Funct
  (\cC^{2},\cC )$.
\end{thm}
We note the following lemma.
\begin{lem} \label{qYB}
Given a category $\cC $ with a premonoidal structure $(\otimes ,\ass
,\com )$, the following quasi--Yang--Baxter diagram commutes.\\ \\
\xymatrix{
\otimes (\otimes \times 1) \ar[r]^{\ass } \ar[d]_{\otimes (\com \times
  1)} & \otimes (1\times \otimes ) \ar[d]^{\otimes (1\times \com )} \\
\otimes (\otimes \times 1)\tau_{(12)} \ar[d]_{\ass \tau_{(12)}} &
\otimes (1\times \otimes )\tau_{(23)} \\
\otimes (1\times \otimes )\tau_{(12)} \ar[d]_{\otimes (1\times \com
  )\tau_{(12)}} & \otimes (\otimes \times 1)\tau_{(23)} \ar[u]_{\ass
    \tau_{(23)}} \ar[d]^{\otimes (\com \times 1)\tau_{(23)}} \\
\otimes (1\times \otimes )\tau_{(23)}\tau_{(12)} & \otimes (\otimes
\times 1)\tau_{(12)}\tau_{(23)} \ar[d]^{\ass \tau_{(12)}\tau_{(23)}} \\
\otimes (\otimes \times 1)\tau_{(23)}\tau_{(12)} \ar[u]^{\ass
  \tau_{(23)}\tau_{(12)}} \ar[d]_{\otimes (\com \times
  1)\tau_{(23)}\tau_{(12)}} & \otimes (1\times \otimes ) \ar[d]^{\otimes 
  (1\times \com )\tau_{(12)}\tau_{(23)}} \\
\otimes (\otimes \times 1)\tau_{(12)}\tau_{(23)}\tau_{(12)} \ar[r]_{\ass 
  } & \otimes (1\times \otimes )\tau_{(23)}\tau_{(12)}\tau_{(23)}}
\end{lem}
{\it Proof} Reading left to right the third row corresponds to $\otimes
(\com \times 1)(\otimes \times 1)$, as does the fourth row. The square
formed is natural. The hexagonal diagrams formed above and below are
those of the definition. Hence the entire diagram commutes.

{\it Proof of Theorem \ref{cohthm4}} The primitive arrows for
interchange (about the region $i$) are of the form $(1,\tau_{i}):(s,\pi
)\rightarrow (s,(i\, i+1)\pi )$ such that $|\vee_{si}s|=2$. In other
words this corresponds to the interchange of two attached (adjacent)
leaves. Every interchange arrow $\tau :(s,\pi )\rightarrow (t,(i\,
i+1)\pi )$ may be written as a sequence of primitive arrows $p_{1}\cdots
p_{m}$ with precisely one corresponding to a primitive interchange
arrow. We define $\Gamma \tau =(\Gamma p_{1})\cdots (\Gamma p_{n})$,
where for an arrow $(f,1):(s,\pi )\rightarrow (t,\pi )$ we have $\Gamma
(s,\pi )=(\Gamma s)\tau_{\pi }$ and $\Gamma (f,1)=(\Gamma f)\tau_{\pi
  }$, and for any primitive interchange $(1,\tau_{i})$
\begin{eqnarray}
\Gamma (1,\tau_{i})=(1^{i-1}\times \com \times 1^{|s|-i+2})\tau_{\pi }
\end{eqnarray}
The proof is completed by showing that this definition is well--defined
and that conditions (\ref{brcond1}) and (\ref{brcond2}) hold whenever
composition is allowed.

We show it is well--defined in two steps. Firstly that there is a
sequence of primitive arrows with precisely one primitive interchange
arrow $(1,\tau_{i}):(s^{\prime },\pi )\rightarrow (t^{\prime },(i\,
i+1)\pi )$ with $s^{\prime -1}(|s|-1)=i$. Secondly, that any two
alternative such sequences form a commutative diagram in $\cC
$. Consider the following diagram.\\ \\
\includegraphics[width=400pt]{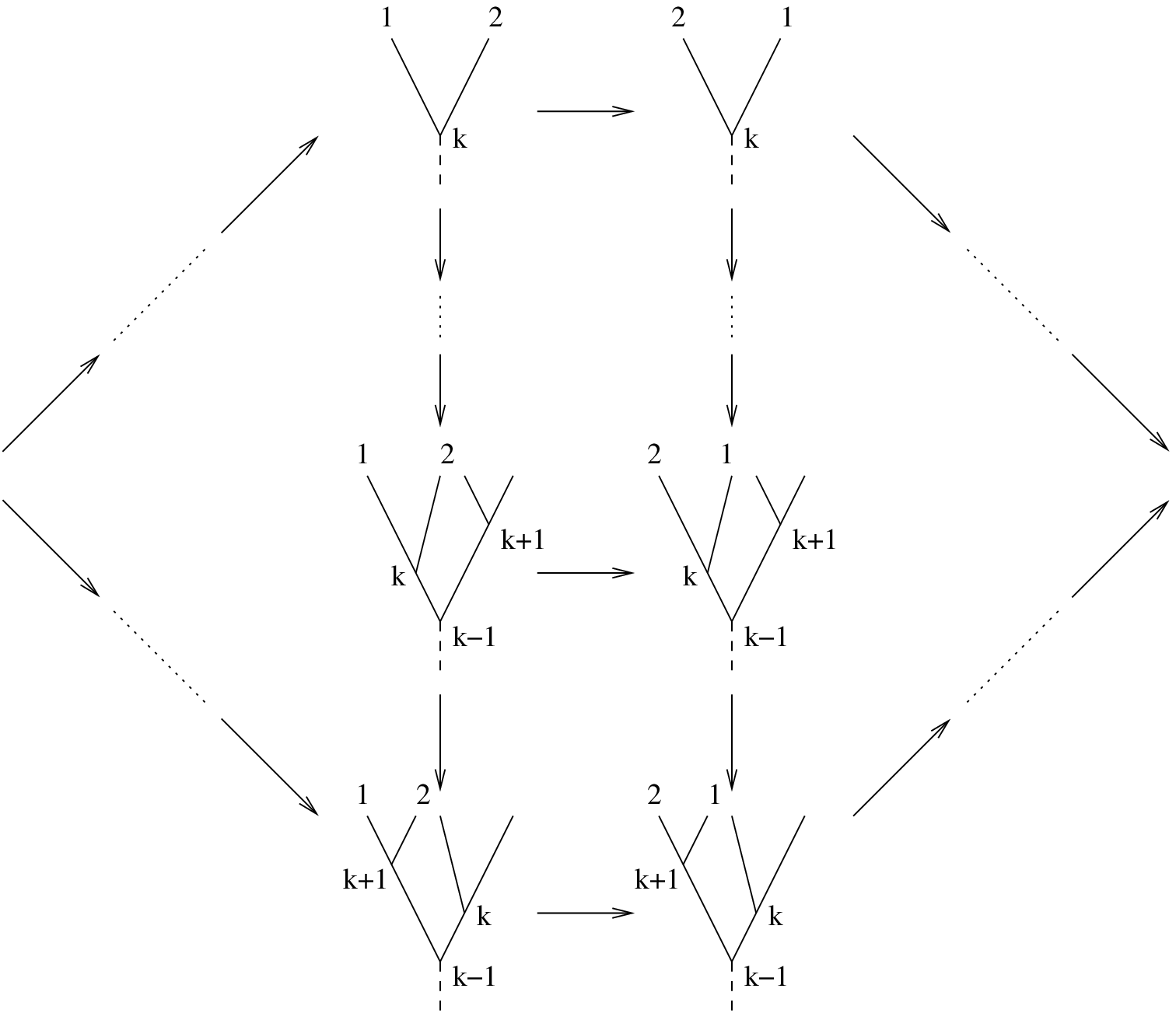}\\ \\
The sequence of arrows running along the top is $p_{1}\rightarrow \cdots
\rightarrow p_{r}$ with the interchange about the $k$th level given by
the centre arrow. We suppose that the level $k+1$ is to the right of $k$
in the source and target trees of this arrow. We construct parallel
sequences of reattachments (vertically downward on diagram) maintaining
the position of the levels greater than $k+1$, into a form containing
the subtree indicated. The diagram is enclosed with the interchange
arrow forming a ladder of natural diagrams under $\Gamma $. Next we
apply the sequence of four arrows corresponding to interchanging the
levels $k$ and $k+1$ completing the region by the relevant diagram of
the definition of braid structure. Finally we can complete the sequence
around the bottom with primitive reattachment arrows. By theorem
\ref{cohthm1} the left and right side diagrams commute. Hence the
sequence around the bottom composes to give the same arrow under $\Gamma
$ as the top sequence. Applying this argument inductively we arrive at a
desired sequence where the interchange occurs about the maximal level
$|s|-1$.

Given an alternative such primitive sequence $q_{1}\rightarrow \cdots
\rightarrow q_{r^{\prime }}$ we can suppose that the primitive
interchange occurs about the $(|s|-1)$th level. Hence we can construct
two parallel sequences of reattachment arrows between the sources and
between the targets of the two interchange arrows preserving the
position of level $|s|-1$. The enclosed diagram commutes under $\Gamma $
because it is a ladder of natural diagrams. Also the remaining two
regions enclosed, one containing the source of $p_{1}$ and $q_{1}$, the
other the target of $p_{r}$ and $q_{r^{\prime }}$, commute under $\Gamma
$ by theorem \ref{cohthm1}. Hence the two sequences give the same arrow
under $\Gamma $ and the definition is well--defined.

Condition (\ref{brcond1}) holds because we can suppose that the source and
target trees of the interchange arrow in a sequence of primitives
composing to give $\tau_{i}$ and $\tau_{j}$ are identical. The result
follows by naturality. Similarly condition (\ref{brcond2}) holds by
lemma \ref{qYB}.

The coherence of the related cases for prebraid and braid
pseudo--monoidal structures are by now a variation on a theme. We make the
following remarks.
\begin{rem}
Theorem \ref{cohthm4} may be weakened to a prebraid structure where the
primitive adjacent interchange arrows are taken as those interchanging two
leaves. The hexagon diagrams define the adjacent interchange of three
leaves. Thus interchanges involving more than three adjacent leaves are
given by iterating the hexagon diagrams.
\end{rem}
\begin{rem}
Alternatively theorem \ref{cohthm4} may be strengthened to the braided
pseudo--monoidal situation where primitive arrows are not restricted by
the requirement that levels are adjacent.
\end{rem}

\section{Braidings and Nodules}                                      %

Finally we bring everything together in the following definition.
\begin{defn}
  A braid premonoidal structure with unit $(\otimes ,\ass ,\com ,\lid
  ,\rid ,e)$ for a category $\cC $ is a premonoidal structure with unit
  $(\otimes ,\ass ,\lid ,\rid ,e)$ and a braid premonoidal structure
  $(\otimes ,\ass ,\com )$.
\end{defn}
Similar definitions hold for prebraid structure with unit and braided
pseudo--monoidal structure with unit.
We define
\begin{eqnarray}
\BNCptr & = & \coprod_{n\in \bW }\Cptr_{n}\times \xB_{n}\times \cN ([n])
\end{eqnarray}
We are now in a position to state the expected coherence result.
\begin{thm} \label{cohthm5}
  If $(\otimes ,\ass ,\com ,\lid ,\rid ,e)$ is a braided premonoidal
  structure with unit for $\cC $ then there is an extension
  of theorems \ref{cohthm3} and \ref{cohthm4} to a functor $\Gamma
  :\BNCptr \rightarrow \Funct (\cC^{\bW },\cC )$.
\end{thm}
The proof of this is very straightforward as are the analogous results
for prebraid and braided pseudo--monoidal structures with unit.

\section{Diagram Calculi}                                            %

Thus far coherence is a functor $\Gamma $ between some groupoid $\Cohr
$, taken as $\Cptr $, $\NCptr $, $\BCptr $ or $\BNCptr $, and $\Funct
(\cC^{\bW },\cC )$. Ultimately, coherence concerns the commutativity of
certain diagrams in $\cC $. Thus we introduce an evaluation functor
\begin{eqnarray}
\ev \equiv \coprod_{n\in \bW }\ev_{n}:\coprod_{n\in \bW }\Funct
(\cC^{n},\cC )\times \cC^{n}\rightarrow \cC
\end{eqnarray}
given by mapping the arrows $(\tau ,f):(F,a)\rightarrow (G,b)$ to
$\ev_{n}(\tau ,f)=(Gf)\tau a$ which by the natural property of $\tau $
is also given by $(\tau b)Ff$. Next we define precisely what we mean by a
diagram in a category.
\begin{defn}
  A collection of arrows $D$ for a category $\cC $ is called a
  (commutative) diagram if given any two composable sequences of arrows
  $f_{1},...,f_{m}$ and $g_{1},...,g_{n}$ from $D$ with matching source
  ($sf_{1}=sg_{1}$) and target ($tf_{m}=tg_{n}$) then we have
  $f_{m}\cdots f_{1}=g_{n}\cdots g_{1}$.
\end{defn}
Clearly if $E\subset D$ and $D$ is a diagram then so is $E$. A diagram
is a labeled directed graph and so inherits the notion of
connectedness. Furthermore, every diagram is the disjoint union of
connected diagrams.
\begin{defn} \label{cohdiag}
  A functor $\Gamma :\cC \rightarrow \clD $ is called coherent
  if for every diagram $D$ of $\clD $ there is a diagram $C$ of
  $\cC $ such that $\Gamma C=D$.
\end{defn}
\begin{rem}
The converse of definition \ref{cohdiag} clearly holds because $\Gamma $ 
is a functor.
\end{rem}
We define the {\it canonical} functor by
\begin{eqnarray}
\can \equiv \coprod_{n\in \bW }\can_{n}:\coprod_{n\in \bW
  }\Cohr_{n}\times \cC^{n}\rightarrow \cC
\end{eqnarray}
where $\can_{n}=\ev_{n}(\Gamma \times 1_{\cC^{n}})$. We can now state
the self evident coherence result.
\begin{thm}
The functor $\can $ is coherent.
\end{thm}

We illustrate the diagram calculus for $\BCptr $. An object $(s,\pi
,\bfa )$ of $\Cptr_{n}\times \xB_{n}\times \cC^{n}$ consists of a
coupling tree $s$, of length $n$ say, with leaves labeled from left to
right by the $n$--tuple of objects $\bfa =(a_{1},...,a_{n})$ from
$\cC^{n}$, and $\pi \in S_{n}$. An arrow $(\sigma ,\tau ,\bff ):(s,\pi
,\bfa )\rightarrow (t,\phi ,\bfb )$ consists of an $n$--tuple of arrows
$\bff =(f_{1},...,f_{n})$ with $f_{i}\in \cC (a_{i},b_{i})$, a
permutation $\sigma \in S_{n-1}$, and a braid $\tau \in \xB_{n}(\pi
,\phi )$. We represent an arrow by a labeled box on a string. Boxes are
free to slide along strings (naturality) and the identity arrow is
simply given by a string. Composition is given by combining vertically
aligned consecutive boxes as depicted in the following diagram for the
composition of
\begin{eqnarray}
((14)(45),\tau_{5}\tau^{-1}_{4}\tau_{3}\tau_{1},f_{1},...,f_{6}):
(01234,1,\bfa )\rightarrow (23014,(23)(56),\bfb )\nonumber
\end{eqnarray}
with
\begin{eqnarray}
((13)(24),\tau^{-1}_{5}\tau_{2},g_{1}, ...,g_{6}):(23014,(23)(56),\bfb
)\rightarrow (14023,(12643),\bfc )\nonumber
\end{eqnarray}
giving
\begin{eqnarray}
(g_{1}f_{1},...,g_{6}f_{6},\tau ,(13452)):(01234,1,\bfa )\rightarrow
(14023,(12643),\bfc )\nonumber
\end{eqnarray}
where $\tau
=\tau_{5}\tau^{-1}_{4}\tau^{-1}_{5}\tau_{3}\tau_{1}\tau_{2}$.\\ \\
\includegraphics[width=400pt]{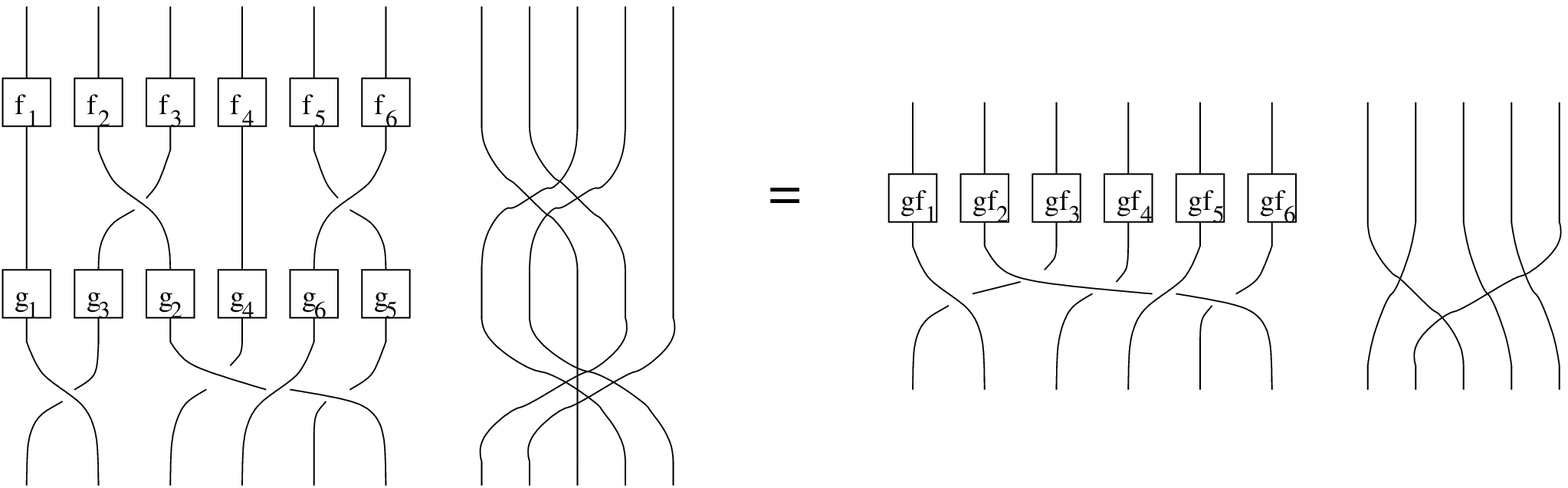}\\ \\


\section*{Acknowledgments}                                          %

WPJ would like to thank Jeremy Martin for reading an early draft. This
research was supported by the New Zealand Foundation for Research,
Science and Technology. Contract number: UOCX0102.


\end{document}